\documentclass[leqno]{amsart}
\usepackage{euscript,amssymb,amsbsy}
\usepackage[arrow,matrix]{xy}

\newtheorem{stuff}{Stuff}[section]
\newtheorem{theorem}[stuff]{\sl Theorem}
\newtheorem{proposition}[stuff]{\sl Proposition}
\newtheorem{lemma}[stuff]{\sl Lemma}
\newtheorem{corollary}[stuff]{\sl Corollary}

\newenvironment{definition}{%
\vskip1ex\refstepcounter{stuff}\trivlist \itemindent 0pt
\item[\hskip\labelsep\sl Definition \thestuff.]%
\ignorespaces}{\endtrivlist\vskip1ex}%
\let\rar\rightarrow
\let\lar\longrightarrow

\let\hra\hookrightarrow
\let\mt\mapsto
\let\lmt\longmapsto
\let\xrar\xrightarrow

\font\tenmsa=msam10 %
\newcommand\hdashpiece{%
{\vrule height2.75pt depth-2.35pt width2.3pt \kern1.7pt}}%
\newcommand\hdashpieces{%
{\hdashpiece\hdashpiece\hdashpiece\hdashpiece}}%
\newcommand\dashar{\mathrel{%
\hdashpieces\kern-0.4pt\hbox{\tenmsa K}}}%

\let\euf\EuScript 
\let\cal\mathcal
\let\mbb\mathbb

\DeclareFontFamily{OT1}{rsfs}{}
\DeclareFontShape{OT1}{rsfs}{n}{it}{<->rsfs10}{}
\DeclareMathAlphabet{\crl}{OT1}{rsfs}{n}{it}

\let\unl\underline

\let\tld\tilde

\let\nit\noindent

\let\disp\displaystyle
\let\srel\stackrel

\let\eps\epsilon
\let\veps\varepsilon

\newcommand\rd{{\rm d}}
\newcommand\Hom{\mathop{\sf Hom}\nolimits}

\newcommand\Mor{\mathop{\sf Mor}\nolimits}
\newcommand\invq{{\slash\kern-0.65ex\slash}}

\numberwithin{equation}{section}

\let\l\lambda
\let\L\Lambda
\let\D\Delta
\let\O\Omega
\let\si\sigma
\let\Si\Sigma
\let\th\theta
\let\sm\setminus
\let\jm\jmath
\newcommand\eL{{\euf L}}

\newcommand\eJ{{\euf J}}
\newcommand\eP{{\euf P}}
\newcommand\ua{{\unl a}}
\newcommand\ud{{\unl d}}
\newcommand\uz{{\unl z}}

\newcommand\Cr{{{({\mbb C}^*)}^r}}
\newcommand\Lie{\mathop{{\cal L}ie}\nolimits}
\newcommand\T{\mathop{\sf T\kern-.2ex}\nolimits}
\newcommand\N{\mathop{\sf N\kern-.1ex}\nolimits}
\newcommand\pr{\mathop{\rm pr}\nolimits}
\begin{document}
\title{Higher genus curves on toric varieties}
\author{Mihai Halic}
\address{Fakult\"at f\"ur Mathematik NA\\ Ruhr Universit\"at Bochum\\ 
44780 Bochum\\ Germany}
\email{halic@ag.ruhr-uni-bochum.de}
\maketitle
\markboth{CURVES ON TORIC VARIETIES}{MIHAI HALIC}

\section*{Introduction}
This article is intended to be an application to the case of 
toric varieties of the ideas developed in my paper \cite{ha}, 
where I have studied the problem of computing the Gromov-Witten 
invariants of quotient varieties.

The quantum cohomology ring of toric varieties was studied in 
\cite{ba}, where is given an explicit formula for the quantum 
multiplication on Fano toric varieties, and also in \cite{sp}, 
where the virtual localization technique is used in order to 
compute genus zero invariants of toric varieties. For the higher 
genus invariants, the combinatorics which appears makes the 
application of the virtual localization formula difficult.

The aim of this article is less ambitious, we just want to have 
an as neat as possible description of the space of morphisms from 
a fixed curve into a smooth and projective toric variety, provided 
the degree of these morphisms is sufficiently large. We impose the 
restriction on the degree because it turns out that in this case 
the space of morphisms in question is smooth and has the expected 
dimension. Eventually, we will compute intersection numbers on a 
certain natural smooth and projective compactification of it. As a 
disclaimer, we should say from the very beginning that we do not 
pretend to describe the stable map compactification \`a la 
Kontsevich-Manin, which is hard to grasp, but instead what we 
find is birational to an irreducible component of this space of 
stable maps.

The article begins with an introductory section whose purpose is 
that of fixing the notations and recalling generalities on toric 
varieties. The study actually starts with the second section which 
treats the problem of compactifying the space of morphisms from a 
smooth and projective curve into a smooth and projective toric 
variety; the conclusions in this direction are contained in 
proposition \ref{inv-quot} and corollary \ref{cor:V}. After the 
description of the generators of the cohomology ring and the 
relations among them, we apply in the last section the 
localization method for computing intersection products on 
our space of morphisms. The general formula given in theorem 
\ref{thm:push-f}, which is the main result of the present article, 
has the shortcoming of being too combinatorial and therefore not 
explicit enough. Neverthless, using it we are able to give in 
proposition \ref{prop:explicit} explicit formulae for certain 
intersection products on one hand and, on the other hand, to 
derive in proposition \ref{prop:0} vanishings induced by the 
primitive collections of the fan defining the toric variety 
we start with.\medskip

\noindent{\it\small Acknowledgment} 
This study was done during my stay at the Maths Department of the 
Ruhr-Universit\"at Bochum, and I take the opportunity to thank 
Prof. H.~Flenner and Prof. A.~Huckleberry for their hospitality.


\section{Setting up the problem}{\label{sct:setup}}

Toric varieties are studied in great detail in \cite{da,fu,oda}, 
but I shall prefer the very nice synthesis of the topic which can 
be found in the third chapter of \cite{ck}. We are starting with 
the notations: $X=X_\Si$ stands for a smooth and projective toric 
variety defined by the fan 
$\Si\subset M^\vee_{\mbb R}:=\Hom_{\mbb Z}(M,{\mbb R})$, with 
$M\cong{\mbb Z}^n$, whose one dimensional faces are denoted 
$\Si(1)$; we let $r:=\#\Si(1)$ and  $l:=r-n$. 
The assumptions on $X$ imply that the integral generators of 
any $n$-dimensional cone of $\Si$ form a ${\mbb Z}$-basis of 
$M^\vee$. We denote $e^1,\dots,e^l,e^{l+1},\dots,e^r$ the 
integral generators of $\Si(1)$, and we are further assuming that 
$(e^{l+1},\dots,e^r)$ actually generate an $n$-dimensional cone. 
Since they form a ${\mbb Z}$-basis of $M^\vee$, there are integers 
$a^\l_\nu$ such that 
\begin{align}{\label{rel-for-gen}}
e^\l +a_\nu^\l e^{l+\nu}=0.
\end{align}
Here we are using the usual summing convention over the indices; in 
the whole paper we let $\l\in\{1,\dots,l\}$, $\nu\in\{1,\dots,n\}$ 
and $\rho\in\{1,\dots,r\}$. We have the exact sequence of tori
\begin{align}{\label{e-sqn}}
1\lar T\srel{\veps}{\lar}{({\mbb C}^*)}^{\Si(1)}\lar S\lar 1,
\end{align}
where $T:=\Hom_{\mbb Z}(A^1(X);{\mbb C}^*)$ and 
$S:=M^\vee\otimes_{\mbb Z}{\mbb C}^*$, which in 
the basis $(e^{l+1},\dots,e^r)$ takes the form 
\begin{align}{\label{exact-sqn}}
1\lar T\srel{\veps}{\lar} {({\mbb C}^*)}^r\lar S\lar 1,
\end{align}
with the homomorphism $\veps$ defined by the characters 
${(\chi_\rho)}_\rho$ as follows:
\begin{align}{\label{eqn:char}}
\chi_\l(t)=t_\l\quad\text{and}\quad
\chi_{l+\nu}(t)
=t^{\ua_\nu}=t_1^{a_\nu^1}\cdot{\dots}\cdot t_l^{a_\nu^l}.
\end{align}
We should keep in mind that this last description depends on 
the choice of a $n$-dimensional cone of $\Si$, and this remark 
will be used over and over in the paper. The homomorphism $\veps$ 
induces a $T$-action on ${\mbb C}^r$, and $X$ is simply the 
quotient for this action. More precisely, there is a $T$-invariant 
open subset $\O\subset{\mbb C}^r$ whose complement 
$Z_X:={\mbb C}^r\sm\O$ has codimension at least two, such that 
$X=\O/T$. In fact $Z_X$ is a union of linear subspaces of 
${\mbb C}^r$,
\begin{align}{\label{Z-X}}
Z_X=\bigcup_\pi{\mbb A}(\pi),
\end{align}
where $\pi$ runs over the set of so-called primitive collections of 
$\Si$ (see sections 1 and 2 of \cite{cox} for definition and proof). 

The generators ${(e^\rho)}_\rho$ of $\Si(1)$ define respectively 
the divisors ${(D_\rho)}_\rho$ on $X$ which, as elements of $A^1(X)$, 
satisfy the linear relations 
\begin{align}{\label{rel-for-div}}
D_{l+\nu}-a_\nu^\l D_\l=0,\text{ for all }\nu=1,\dots,n.
\end{align}
In this way we get an isomorphism 
$A^1(X)\cong {\mbb Z}D_1\oplus\dots\oplus{\mbb Z}D_l$, 
determined by the choice of the $n$-dimensional cone 
$(e^{l+1},\dots,e^r)$ of $\Si$.

\begin{definition}{\label{def:deg}}
For $C$ a smooth and projective curve of genus $g$, we say that a 
morphism $u:C\rar X$ has multi-degree $\ud={(d_\rho)}_\rho$ if 
$d_\rho=D_\rho\cdot u_*C$ for $\rho=1,\dots,r$.
\end{definition}

Of course, in definition above the integers $d_\rho$ are not 
independent but are related by $d_{l+\nu}=a_\nu^\l d_\l$ for 
all $\nu$. 

The following well-known `Euler sequence' on $X$, 
\begin{align}{\label{euler}}
0\lar {\cal O}_X^{\oplus l}\lar 
\bigoplus_\rho{\cal O}(D_\rho)\lar
\T_X\lar 0,
\end{align}
and the Riemann-Roch theorem immediately implies the 

\begin{lemma}{\label{lm:smooth-map}}
The space $\Mor_\ud(C,X)$ of morphisms from $C$ to $X$ having 
multi-degree $\ud$ is smooth as soon as $d_\rho>2g-1$ for all 
$\rho=1,\dots,r$. In this case, it is also irreducible and has 
the expected dimension 
\begin{align}{\label{exp-dim}}
\dim\Mor_\ud(C,X)=\sum_\rho d_\rho - n(g-1).
\end{align}
\end{lemma}

\begin{proof}
All the statements are obvious except the one concerning the 
irreducibility of $\Mor_\ud(C,X)$ which will be proved in 
corollary \ref{cor:V}.
\end{proof}

As this first section is devoted to set up the stage, let us 
mention a probably well-known generality about toric varieties 
which will be needed in the proof of proposition \ref{inv-quot}.

\begin{proposition}{\label{inv-t}}
Let $X$ be a smooth and projective toric variety as before. There 
is an ample line bundle $A\rar{\mbb C}^r$ having the properties:

{\rm (i)\,} it linearizes the standard $\Cr$-action on ${\mbb C}^r$;

{\rm (i)\,} the corresponding set of $T$-semi-stable points is 
precisely $\O\subset{\mbb C}^r$.
\end{proposition}

\begin{proof} Since $X$ is projective, there are characters 
$\beta$ of $T$ such that the associated line bundles 
$A_\beta:=\O\times_\beta{\mbb C}\rar X$ are very ample. Sections 
of $A_\beta$ can be naturally identified with $T$-equivariant 
holomorphic functions on $\O$, and {\it a fortiori} on ${\mbb C}^r$ 
because $Z_X={\mbb C}^r\sm\O$ has codimension at least two. We say 
that $f:{\mbb C}^r\rar{\mbb C}$ is $T$-equivariant if 
$$
f\bigl(\veps(t')\times z\bigr)=\beta(t')f(z),\,
\forall t'\in T\text{ and }\forall z\in{\mbb C}^r,
$$
and denote ${\cal O}({\mbb C}^r)_\beta$ the (finite dimensional) 
vector space of such functions. There is a natural map 
$$
\begin{array}{l}
F:{\mbb C}^r\lar {{\cal O}({\mbb C}^r)_\beta}^\vee,\\
\kern2em z\lmt F(z):\kern2ex 
\langle F(z),f\rangle=f(z),
\end{array}
$$
having the properties:

$\bullet$ it covers the projective embedding of $X$ defined by 
the linear system of $A_\beta$;

$\bullet$ $F(z)=0$ for all $z\in Z_X$ and 
$F(z)\neq 0$ for all $z\in\O$;

$\bullet$ $F\bigl(\veps(t')\times z\bigr)=\beta(t')F(z)$ 
for all $t'\in T$ and $z\in{\mbb C}^r$.

\nit In a certain basis of ${{\cal O}({\mbb C}^r)_\beta}^\vee$, 
the $\Cr$-action can be diagonalized, and we denote $\mu_1,\dots,
\mu_\omega$ the corresponding characters; they have the property 
that $\mu_k\circ\veps=\beta$ for $k=1,\dots,\omega$. With respect 
to this basis, we are defining the $\Cr$-action on ${\mbb P}\bigl(
{{\cal O}({\mbb C}^r)_\beta}^\vee\oplus{\mbb C}\bigr)$, together 
with a linearization in the standard ample line bundle over this 
projective space as follows: $\Cr$ acts on 
${{\cal O}({\mbb C}^r)_\beta}^\vee$ by the 
characters $1,\mu_1^{-1}\mu_2,\dots,\mu_1^{-1}\mu_\omega$ 
and on the ${\mbb C}$-term by $\mu_1^{-1}$. 
Similarly, we let $\Cr$ act on ${\mbb C}^{r+1}={\mbb C}^r
\oplus{\mbb C}$, in the standard fashion on ${\mbb C}^r$ 
and trivially on the last component. 

The ample line bundle $A\rar{\mbb C}^r$ we are looking for 
is obtained by restricting via 
$$
{\mbb C}^r\cong{\rm Graph}(F)\subset
{\mbb P}\bigl({\mbb C}^r\oplus{\mbb C}\bigr)
\times
{\mbb P}\bigl(
{{\cal O}({\mbb C}^r)_\beta}^\vee\oplus{\mbb C}
\bigr)
$$
the natural (tensor product) ample line bundle on the product 
of these projective spaces. After replacing $\beta$ with 
a sufficiently large power of it, one may check that the 
points $(z,1)\times (0,1)$, $z\in {\mbb C}^r$, are $T$-unstable, 
while the points $(z,1)\times (y,1)$, $y\neq 0$, are 
$T$-semi-stable. Consequently the $T$-semi-stable set of 
${\mbb C}^r$ is precisely $\O$. 
\end{proof}

We continue the preparatory material by fixing once for 
all a point $\zeta_0\in C$ and we consider the Poincar\'e 
bundle $\eL_0\rar\eJ\times C$ whose restriction 
${\eL_0}|_{\eJ\times\{\zeta_0\}}={\cal O}_\eJ$; for $d\in{\mbb Z}$, 
we define $\eL_d:=\eL_0\otimes\pr_C^*{\cal O}(d\zeta_0)$. 
As one expects, $\eJ$ denotes the Jacobian variety of $C$ and, for 
all integers $d$, $\eL_d$ parameterizes line bundles of degree $d$ 
over the curve $C$. 

The topological type of a holomorphic principal $T$-bundle over 
$C$ is determined by its multi-degree $\ud'=(d_1,\dots,d_l)$; 
holomorphic principal $T$-bundles over $C$ with fixed multi-degree, 
are parameterized by the $l^{\rm th}$ power of the Jacobian of $C$. 
We denote 
$$
\eP_{\ud'}\lar \eJ^l\times C
$$
the universal principal $T$-bundle parameterizing principal 
$T$-bundles over $C$ with multi-degree $\ud'$, trivialized at 
$\zeta_0$. The bundle $\eP_{\ud'}$ is uniquely determined by 
the choice of the Poincar\'e bundles $\eL_d$ above.

Later on we will see that is useful to consider the principal bundle 
$$
\eP_\ud\lar\eJ^r\times C
$$
which parameterizes ${({\mbb C}^*)}^r$-bundles over 
$C$ having multi-degree $\ud=(d_1,\dots,d_r)$. For 
\begin{align}{\label{map-phi}}
\begin{array}{l}
\psi:\eJ^l\lar\eJ^r\quad\text{defined by}\\[1.5ex] 
(L_1,\dots,L_l)\lmt 
(L_1,\dots,L_l,a_1^\l L_\l,\dots,a_n^\l L_\l),
\end{array}
\end{align}
we get the commutative diagram
\begin{align}{\label{cd1}}
\xymatrix{
\eP_{\ud'}\times_\veps{({\mbb C}^*)}^r=\psi^*\eP_\ud
\ar[d]\ar[r]& \eP_\ud\ar[d]\\ 
\eJ^l\times C\ar[r]^-\psi &\eJ^r\times C.
}
\end{align}
Shortly, the reason for introducing this new ingredient is that 
for writing the left-hand-side of \eqref{cd1} we have chosen a 
cone of $\Si$, while the right-hand-side is symmetric, the 
information on the structure of $\Si$ being encoded in the map 
$\psi$.


\section{Description of the space of morphisms}{\label{sct:descr}}

All the subsequent constructions are motivated by the following 
very simple remark: given a morphism $u:C\rar X$ having multi-degree 
$\ud$, the pull-back $P:=u^*\Omega\rar C$ is a holomorphic principal 
$T$-bundle whose multi-degree is $\ud'$ (this is easy to see). 
The morphism $C=P/T\rar P\times_T{\mbb C}^r$ is just a section 
of a rank $r$ vector bundle over $C$ in which the torus $T$ 
still acts, covering the identity of $C$. Any two sections which 
are in the same $T$-orbit give rise to the same morphism from $C$ 
into $X$ (some care is actually required at this point). This is 
the sort of correspondence which will be exploited in this 
section. The inequalities $d_\rho>2g-1$ appearing in lemma 
\ref{lm:smooth-map} will be assumed in the rest of the paper. 

We start with the vector bundles associated to \eqref{cd1}
\begin{align}{\label{cd2}}
\xymatrix{
\eP_{\ud'}\times_T{\mbb C}^r
=\psi^*(\eP_\ud\times_\Cr {\mbb C}^r)
\ar[d] \ar[r] & 
\eP_\ud\times_\Cr {\mbb C}^r
=\oplus_\rho \eL_\rho \ar[d]\\ 
\eJ^l\times C\ar[r]^-\psi & \eJ^r\times C,
}
\end{align}
and notice that $\psi^*\eL_\rho=\eP_{\ud'}\times_{\chi_\rho}
{\mbb C}=:\crl L_\rho$ (the characters $\chi_\rho$ are defined 
by \eqref{eqn:char}). Taking the direct images 
\begin{align}{\label{cd3}}
\xymatrix{
\crl V=\oplus_\rho\crl V_\rho
:=p_*(\oplus_\rho\crl L_\rho)=\psi^*\crl W
\ar[d]\ar[r] &
\crl W=\oplus_\rho\crl W_\rho
:=p_*(\oplus_\rho \eL_\rho)
\ar[d]\\ 
\eJ^l\ar[r]^-\psi& \;\eJ^r,
}
\end{align}
we recognize in $\crl W_\rho$ the Picard vector bundles 
associated respectively to the Poincar\'e bundles $\eL_\rho$; 
the rank of $\crl W$ is given by the formula 
$$
{\rm rk\,}\crl W=\sum_\rho d_\rho-r(g-1).
$$
The action of $T$ on ${\mbb C}^r$ induces actions on 
$\oplus_\rho\crl L_\rho$ and $\oplus_\rho\eL_\rho$ covering 
respectively the identities of $\eJ^l\times C$ and $\eJ^r\times C$ 
and, {\it a fortiori}, there are natural $T$-actions on $\crl V$ 
and $\crl W$ which cover respectively the identities of $\eJ^l$ 
and $\eJ^r$ and moreover preserve the decompositions 
$\crl V=\oplus_\rho\crl V_\rho$ and 
$\crl W=\oplus_\rho\crl W_\rho$. 

The remark at the beginning of this section tells that 
the space of morphisms from $C$ to $X$ should be the quotient 
`$\crl V/T$'. Of course, this should not be taken {\it ad litteram} 
but in the spirit of geometric invariant theory. What we shall 
actually construct is the quotient of $\crl W$ for the $T$-action, 
and $\crl V\invq T$ will be just its pull-back by $\psi$. 

One can spot at the first glance a `nice' Zariski open subset of 
$\crl W$ on which $T$ acts freely
\begin{align}{\label{stable-pts}}
\crl W^o:=\{ s\in\crl W \mid 
{\rm Image\,}s\not\subset\eP_\ud\times_\Cr Z_X\}.
\end{align}
As the closed subvariety $Z_X\subset{\mbb C}^r$ which had to be 
`thrown away' for obtaining $X$ was a union of coordinate subspaces, 
\begin{align}{\label{Z-W}}
Z_W:=
\Gamma\bigl(
C,\eP_\ud\times_\Cr Z_X
\bigr)
=\bigcup_\pi\Gamma\bigl(
C,\eP_\ud\times_\Cr{\mbb A}(\pi)
\bigr)
\end{align}
is still a union of subvector bundles of $\crl W$, and $\crl W^o
=\crl W\sm Z_W$. Even if $T$ is acting freely on $\crl W^o$, it is 
possibly not so clear that the quotient $\crl W^o/T$ exists as a 
complex manifold. 

\begin{lemma}{\label{lm:quot}}
$\crl W^o/T$ has a natural structure of a Hausdorff complex 
analytic variety.
\end{lemma}

\begin{proof}
As expected, we put on $\crl W^o/T$ the final topology for 
the projection $\crl W^o\rar\crl W^o/T$; we must prove that 
the quotient is Hausdorff when we consider on $\crl W^o$ the 
analytic topology given by small balls. For two sections 
$s,s'\in\crl W^o$ lying over the same point in $\eJ^r$, which are 
not in the same $T$-orbit, we want to prove that there are 
neighborhoods ${\cal U}\ni s$ and ${\cal U}'\ni s'$ such 
that $T{\cal U}\cap T{\cal U}'=\emptyset$ (when $s$ and $s'$ 
lie above two different points of $\eJ^r$ everything is clear). 
From the very definition of $\crl W^o$ we deduce that 
$$
s_\zeta,s'_\zeta\in \Omega,\quad\forall\zeta\in C^o,
$$
with $C\sm C^o$ a finite set. We are distinguishing two cases: 

Case (1)\; When there is a point $\zeta\in C^o$ such 
that $s_\zeta$ and $s'_\zeta$ are not in the same $T$-orbit 
in $\O$, the existence of the two disjoint $T$-invariant 
neighborhoods of $s$ and $s'$ is immediate;

Case (2)\; It might happen that for any $\zeta\in C^o$, the 
evaluations $s_\zeta$ and $s'_\zeta$ are in the same $T$-orbit 
(such a situation does appear in the simple case when we projectivize 
the space of sections of a line bundle). In this case, since $T$ acts 
freely on $\Omega$, there is a morphism $\tau:C^o\rar T$ such that 
$s'_\zeta=\tau_\zeta\cdot s_\zeta$. If there are no neighborhoods 
${\cal U}$ and ${\cal U}'$ as wanted, we deduce the existence of 
sequences ${(s_k)}_k\subset\crl W^o$ and ${(t_k)}_k\subset T$ such 
that $s_k\xrar{\| \cdot \| }s$ and 
$t_k\cdot s_k\xrar{\| \cdot \| }s'$, where the norm 
$\| \cdot \|$ on $\crl W$ is defined by 
$$
\|s\|:=\max_{\zeta\in C}\|s_\zeta\| .
$$
Let us fix a closed disk $\bar\D\in C^o$. For positive $\eps$, 
there is a rank $k_\eps$ such that for $k\geq k_\eps$, 
$$
s_{k,\zeta}\in\O,\;\forall\zeta\in\bar\D
\quad\text{and}\quad
\| s'_\zeta - t_k\cdot s_k\| < \eps.
$$ 
It follows that 
$| \tau_\zeta\cdot s_\zeta - t_k\cdot s_{k,\zeta} | <\eps,
\;\forall\zeta\in\bar\D$, which in turn implies that 
$| s_\zeta - \tau_\zeta^{-1}t_k\cdot s_{k,\zeta}|<\eps$ 
for a possibly different choice of $\eps$ (this is because we 
have restricted ourselves to the compact $\bar\D$). Since for 
$\zeta\in\bar\D$ the $s_k$'s take values in $\O$, we deduce that 
$t_k\xrar{k\rar\infty}\tau_\zeta$ for all $\zeta\in\bar\D$; 
but this means that $\tau:C^o\rar T$ is constant when 
restricted to $\D$ and consequently is constant everywhere. 
We conclude that $s$ and $s'$ are in the same $T$-orbit 
in $\crl W$, which contradicts our assumption.
\end{proof}

This lemma solves the problem of the $T$-action on $\crl W^o$, but 
one would like to work with compact manifolds, and no one guarantees 
that $\crl W^o/T$ is so. The purpose of the next proposition is to 
prove that $\crl W^o/T$ is actually a projective, and not just an 
analytic variety. 

\begin{proposition}{\label{inv-quot}}
The $T$-action on $\crl W$ can be linearized in an ample 
line bundle over it, such that the semi-stable set for this 
action coincides with $\crl W^o$ defined by \eqref{stable-pts}. 
Moreover, the invariant quotient $W:=\crl W\invq T$ is a smooth 
and projective variety of dimension 
$$
\dim W=\sum_\rho d_\rho +n.
$$
\end{proposition}

\begin{proof}
According to proposition \ref{inv-t} there is an ample line bundle 
on $A\rar{\mbb C}^r$ which linearizes the $\Cr$-action on 
${\mbb C}^r$, such that the corresponding $T$-semi-stable set is 
$\O$, and consequently $X={\mbb C}^r\invq T$. This one induces the 
relatively ample line bundle 
$$
\bar A:=\eP_\ud\times_\Cr A\rar\eP_\ud\times_\Cr{\mbb C}^r,
$$ 
and the action of $T$ is still linearized in $\bar A$. Tensoring 
$\bar A$ with a sufficiently ample line bundle on $\eJ^r\times C$ 
we obtain an ample line bundle on $\eP_\ud\times_\Cr {\mbb C}^r$, 
together with a linearization of the $T$-action in it, which has 
the additional property that the $T$-semi-stable locus is precisely 
$\eP_\ud\times_\Cr \Omega$. Identifying 
$\crl W=p_*\bigl(\eP_\ud\times_\Cr {\mbb C}^r\bigr)$ with the space of 
morphisms from $C$ into $\eP_\ud\times_\Cr {\mbb C}^r$ which represent 
the class of a section, we find ourselves in the situation studied in 
section 2 of \cite{ha}, where is shown that in this case is possible 
to linearize (rather canonically) the $T$-action on $\crl W$ in an 
ample line bundle such that the corresponding semi-stable points 
have the property that their image is not completely contained 
in the unstable locus of $\eP_\ud\times_\Cr {\mbb C}^r$ (according to 
corollary 2.4 {\it loc.cit.}). Denoting $\crl W^s$ the set of 
$T$-semi-stable points of $\crl W$, we have found that 
$\crl W^s\subset\crl W^o$.

We want to prove now that $W=\crl W^s/T$ is projective; 
quasi-projectivity comes for free from the very construction, 
so that remains to prove the compactness (completeness). Since 
the $T$-action on $\crl W$ covers the identity of $\eJ^r$, 
$\crl W^s/T$ comes with the projection 
$$
q:\crl W^s/T\lar \eJ^r.
$$
The first claim is that this map is surjective: indeed, according 
to theorem 2.5 in \cite{ha}, points $s\in\crl W$ whose image is 
contained in $\eP_\ud\times_\Cr\O$ are $T$-semi-stable. Since the 
`bad' set $Z_X\subset{\mbb C}^r$ has codimension at least two and 
because the line bundles $\eL_\rho$ are globally generated, we 
deduce that for any $j\in\eJ^r$ there is a section 
$s\in \Gamma(C,{\eP_\ud\times_\Cr {\mbb C}^r}|_j)$ having 
the property that its image is disjoint from $\eP\times_\Cr Z_X$. 
This proves the claim. 

Consequently is enough to prove the compactness of the fibres of the 
projection $q$. The following lemma describes the fibrewise situation.  

\begin{lemma}{\label{lm:comp}}
Suppose we are given a torus action on ${\mbb C}^r$ as described 
at the beginning of section \ref{sct:setup}, so that the quotient 
is a smooth and projective toric variety. Let us consider now the 
action $T\times{\mbb C}^R\rar{\mbb C}^R$, where
$$
{\mbb C}^R
:={\mbb C}^{N_1}\oplus\dots\oplus{\mbb C}^{N_r}
$$
and the torus $T$ acts on the direct summands of ${\mbb C}^R$ 
respectively by the same characters ${(\chi_\rho)}_\rho$  as on 
${\mbb C}^r$. Then the quotient is a smooth, compact toric variety. 
\end{lemma}

\begin{proof}
The compactness of the quotient can be easily seen using the moment 
map description of toric varieties. In coordinates, the moment map 
corresponding to the $T$-action on ${\mbb C}^r$ is (see section 3.3 
in \cite{ck}):
\begin{align}{\label{m-map}}
m:{\mbb C}^{n+l}\lar{\mbb R}^l,\quad 
m(z)=\frac{1}{2}\left(
\begin{array}{c}
{|z^1|}^2+a_\nu^1{|z^{l+\nu}|}^2\\ 
\vdots\\ 
{|z^l|}^2+a_\nu^l{|z^{l+\nu}|}^2
\end{array}
\right).
\end{align}
Then $X$ can be described as $X=m^{-1}(a)/T_{\mbb R}$ for (a well 
chosen) $a\in {\mbb R}^l$, where $T_{\mbb R}=(S^1)^l$ denotes the 
real torus. Since $X$ is compact, $m^{-1}(a)$ is still compact and 
$T_{\mbb R}$ acts freely on it. 

Let us move now to the new situation and denote $\uz^\rho$  the 
points of ${\mbb C}^{N_\rho}$. The moment map in this case has 
the form 
\begin{align}{\label{M-map}}
\crl M:{\mbb C}^R\lar{\mbb R}^l,\quad 
\crl M(\uz)=\frac{1}{2}\left(
\begin{array}{c}
{|\uz^1|}^2+a_\nu^1{|\uz^{l+\nu}|}^2\\ 
\vdots\\ 
{|\uz^l|}^2+a_\nu^l{|\uz^{l+\nu}|}^2
\end{array}
\right)
\end{align}
and we want to prove that the quotient 
$Y:=\crl M^{-1}(a)/T_{\mbb R}$ is smooth and compact. 
It is easy to see that the maps 
$$
m^{-1}(a)\ni (z^1,\dots,z^r)\lmt 
(|z^1|,\dots,|z^r|)\in {\mbb R}^r
$$
and 
$$
\crl M^{-1}(a)\ni (\uz^1,\dots,\uz^r)\lmt 
(|\uz^1|,\dots,|\uz^r|)\in {\mbb R}^r.
$$
have the same image, which is compact since $m^{-1}(a)$ is so. The 
compactness of $\crl M^{-1}(a)$ is implied now by the compactness 
of standard spheres. The action of $T_{\mbb R}$ on $\crl M^{-1}(a)$ 
is free because for any point $\uz_o$ which solves the equation 
$\crl M(\uz)=a$, is possible to find a $({\mbb C}^*)^r$-equivariant 
embedding ${\mbb C}^r\hra {\mbb C}^R$ which pass through  $\uz_o$. 
The conclusion follows now from the fact that $T_{\mbb R}$ acts 
freely on $m^{-1}(a)$. 
\end{proof}

Coming back to our proposition, we deduce from the lemma that the 
fibres of $q$ are smooth and compact toric varieties, all isomorphic 
to $Y=\crl M^{-1}(a)/T_{\mbb R}$. Since $\crl W^s/T$ is 
quasi-projective, it follows that $Y$ is actually projective and so 
is the invariant quotient $W=\crl W^s/T$. 

Remains to prove that $\crl W^s=\crl W^o$: the natural inclusion 
$\crl W^s/T\hra \crl W^o/T$ being an open map, its image is both 
open and closed in $\crl W^o/T$, so that $\crl W^s/T =\crl W^o/T$ 
and therefore $\crl W^s=\crl W^o$.
\end{proof}

According to \eqref{M-map}, the moment map induced by the 
$T_{\mbb R}$-action on $\crl W$ is given by 
\begin{align}{\label{moment-map}}
\crl M(s)=\frac{1}{2}\left(
\begin{array}{c}
\disp\int_C \bigl(
{|s^1_\zeta|}^2+a_\nu^1{|s^{l+\nu}_\zeta|}^2
\bigr)\rd\gamma(\zeta)\\ 
\vdots\\ 
\disp\int_C 
\bigl( 
{|s^l_\zeta|}^2+a_\nu^l{|s^{l+\nu}_\zeta|}^2
\bigr)\rd\gamma(\zeta)
\end{array}
\right)=\int_C m(s_\zeta)\rd\gamma(\zeta),
\end{align}
where $s={(s^\rho)}_\rho\in\crl W$ and $\rd\gamma$ denotes a volume 
form on $C$.  This formula is in agreement with the computations 
done in section 3 of \cite{ha}, namely with equation (3.3) in 
there. Actually, this is the reason why we consider the level set 
$\{\crl M=a\}\subset\crl W$, and not another one, for describing 
the invariant quotient $\crl W\invq T$.

We recall now that what we are actually interested in is a 
compactification of the space $\Mor_\ud(C,X)$. The reason for 
introducing the variety $W$ was to have a `symmetric object' in our 
hands, in the sense that it does not depend on the choice of some 
particular cone of $\Si$. The compactification we are looking for is 
$V:=\crl V\invq T$ (see \eqref{cd3}), which can be now easily described 
as $V=\psi^*W=\eJ^l\times_{\eJ^r}W$. It is a fibre space over $\eJ^l$, 
with all the fibres isomorphic to the toric variety $Y$ constructed 
in lemma \ref{lm:comp}. We collect this information in the 

\begin{corollary}{\label{cor:V}}
The space of morphisms $\Mor_\ud(C,X)$ is irreducible and the variety 
$V:=\psi^*W$ is a smooth and projective compactification of it.
\end{corollary}

Now becomes clear our statement in the introduction, that $V$ is 
definitely not the stable map compactification of the space of 
morphisms from $C$ to $X$, but is only birational to an irreducible 
component of this later. Indeed, the space of stable maps whose 
stabilized domain is $C$ contains, when $g\geq 2$, the component 
whose points correspond to the following morphisms: the domain of 
definition is the singular curve consisting of $C$ with ${\mbb P}^1$ 
attached at some point; the map is constant on $C$ and has 
multi-degree $\ud$ on ${\mbb P}^1$. Is also true that this 
component has strictly larger dimension than the expected 
one, and therefore is not clear how does it contribute to 
the Gromov-Witten invariants.


\section{Cohomology of the space of morphisms}{\label{sct:cohom}}

Since the projection $q:W\rar\eJ^r$ is a fibre bundle, the 
Leray-Hirsch theorem says that the cohomology of $W$ is generated by 
the the cohomology of $\eJ^r$ and the cohomology of the fibre $Y$ (a 
similar remark applies to $q:V\rar\eJ^l$). Let us define now the line 
bundles
\begin{align}{\label{line-bdl}}
\L_\rho:=\crl W^o\times_{\chi_\rho}{\mbb C}\lar W,\quad
\forall\rho=1,\dots,r.
\end{align}
on $W$. The interest comes from the fact that there is a rational 
evaluation map 
\begin{align}{\label{ev}}
ev:V\times C\dashar X,
\end{align}
and the $\psi^*\L_\rho$'s coincide respectively with the pull-backs 
under $ev$ of the line bundles ${\cal O}(D_\rho)$ on $X$, at least 
on the domain of definition of $ev$. Since the classes $D_\rho$ 
generate the cohomology of $X$, we may hope that integrals as 
$\int_V\prod_\rho{(\psi^*\L_\rho)}^{m_\rho}$ are related to 
enumerative invariants of $X$. Morally, they should count the 
number of morphisms from $C$ to $X$ satisfying certain incidence 
conditions.

\begin{proposition}{\label{prop:gen-cohom}}
The integral cohomology of $W$ is generated as a ring by the 
integral cohomology of $\eJ^r$ and the classes $\L_1,\dots,\L_r$.
\end{proposition}

\begin{proof}
The statement follows from the fact that when restricted to 
the fibres of $q$, the classes $\L_1,\dots,\L_r$  generate the 
integral cohomology ring of $Y$ (recall that $Y$ is a toric variety). 
\end{proof}

For making computations, we must find the relations among the classes 
$\L_\rho$. The obvious relations are the linear ones like
\begin{align}{\label{lin-rel}}
\L_{l+\nu}=a_\nu^\l\L_\l,\quad\forall\nu\in\{1,\dots,n\}.
\end{align}
and the others corresponding to the remaining $n$-dimensional cones of 
$\Si$. In analogy with the case of toric varieties, we are going to 
describe the non-linear relations among the $\L_\rho$'s which arise 
from the primitive collections of $\Si$. We start noticing that the 
`bad' sets which must be removed for constructing the quotients $X$ 
and $W$ behave in a rather functorial way: indeed, according to  
\eqref{Z-W}
$$
Z_W=\Gamma\bigl(C,\eP_\ud\times_\Cr Z_X)
=\bigcup_\pi\Gamma\bigl(C,\eP_\ud\times_\Cr{\mbb A}(\pi)\bigr),
$$ 
which is a union of linear subvector bundles of $\crl W$; as usual, 
$\pi\subset\{1,\dots,r\}$ runs over the primitive collections of 
the fan defining $X$. This means that the primitive collections 
of $W$, more precisely the primitive collections of the fan defining 
the toric fibre $Y$ of $q:W\rar\eJ^r$, are simply $\ker(\pr_\pi)$, 
for $\pr_\pi:\crl W\lar \oplus_{\rho\in\pi}\crl W_\rho$. We deduce 
that for any primitive collection $\pi$ of the fan defining $X$, 
\begin{align}{\label{rel-Y}}
\jmath^*_Y
\left(
\prod_{\rho\in\pi}\L_\rho^{N_\rho}
\right)
=0,\text{ for }N_\rho:={\rm rk\,}\crl W_\rho=d_\rho-(g-1).
\end{align}

However, we would like to have expressions for the products 
$\prod_{\rho\in\pi}\L_\rho^{N_\rho}$ as elements of 
$H^*(W)\cong H^*(\eJ^r)\otimes H^*(Y)$. To begin with, we observe 
that for each $\rho$ there is a {\it sheaf} monomorphism 
$$
0\rar \L_\rho^{-1}\rar q^*\crl W_\rho
\quad\text{given by}\quad 
[s,z]\lmt z\pr_\rho s,\;
\forall s\in\crl W^o\text{ and }\forall z\in{\mbb C}.
$$
Equivalently, one can say that for each $\rho$ there is a canonical 
non-zero section $0\rar{\cal O}_W\rar q^*\crl W_\rho\otimes\L_\rho$. 
From the description of $Z_W$ it follows that for every primitive 
collection $\pi$,
\begin{align}{\label{pi-euler}}
0\lar{\cal O}_W\lar\bigoplus_{\rho\in\pi}
q^*\crl W_\rho\otimes\L_\rho
\end{align}
is a monomorphism of vector bundles, and consequently

\begin{proposition}{\label{prop:rel}}
For any primitive collection $\pi$ of the fan defining $X$, the 
Euler class
$$
e
\bigl(
\oplus_{\rho\in\pi} 
q^*\crl W_\rho\otimes\L_\rho
\bigr)
=0.
$$
\end{proposition}

The vanishing \eqref{rel-Y} is an immediate consequence, and when 
$X$ is a projective space, so that $W={\mbb P}(\crl W)$, proposition 
\ref{prop:rel} reduces to the standard Grothendieck relation for 
${\cal O}_{{\mbb P}(\crl W)}(1)\rar W$. It is probably true 
that the linear relations \eqref{lin-rel} and the non-linear 
ones in proposition \ref{prop:rel} generate the ideal of all 
the relations among the $\L_\rho$'s.


\section{Localization}{\label{sct:local}}

In this section we will apply the localization method developed in 
\cite{ab} for computing intersection products of the $\L_\rho$'s 
on $W$, the ultimate goal being to compute intersection numbers 
on $V$. We will apply the localization formula in cohomology with 
respect to the action of the torus $S$ defined by the exact sequence 
\eqref{e-sqn}. 

First of all we have to make explicit the $S$-action on $W$ and to 
describe the corresponding fixed point set $W^S$. Since $Z_W$ defined 
by \eqref{Z-W} is $\Cr$-invariant, its complement $\crl W^o$ is still 
$\Cr$-invariant and therefore $S=\Cr/T$ acts on $W=\crl W^o/T$ and 
moreover the evaluation map \eqref{ev} is $S$-equivariant. This 
remark implies that if $[s]\in W^S$, then for all $\zeta\in C$ such 
that $s(\zeta)\not\in Z_W$, $ev_{[s]}(\zeta)\in X^S$. But $X^S$ 
consists of finitely many points: they correspond in a bijective 
fashion to the $n$-dimensional cones of $\Si$ and their number equals 
the Euler characteristic of $X$. For $x\in X^S$, we shall denote 
$\si_n(x)$ the corresponding $n$-dimensional cone of $\Si$, and by 
$\bar O_x\subset{\mbb C}^r$ the closure of the $T$-orbit above $x$. 
In fact $\bar O_x$ is the linear $l$-dimensional subspace of 
${\mbb C}^r$ defined by the equations
$$
\bar O_x=\{ z^\rho=0\mid \rho\in\si_n(x)\}.
$$
Our discussion implies that for any $[s]\in W^S$ the image of the 
evaluation $ev_{[s]}$ is a point $x\in X^S$, and this in turn means 
that $s\in\Gamma(C,\eP_\ud\times_\Cr\bar O_x)$. What we have obtained 
so far is that 
$$
W^S\subset
\bigcup_{x\in X^S}\Gamma (C,\eP_\ud\times_\Cr\bar O_x)^o/T
=:\bigcup_{x\in X^S}W(x),
$$
and our goal is to show that this inclusion is in fact an equality. 

We are going to check that the component 
$\Gamma (C,\eP_\ud\times_\Cr\bar O_{x_0})^o/T$ is fixed by $S$, 
for $x_0\in X^S$ the point corresponding to the cone 
$\si_n(x_0)=\langle e_{l+1},\dots,e_r\rangle$ of $\Si$. With 
respect to this choice of coordinates, the $T$-action on 
${\mbb C}^r$ is given by \eqref{eqn:char} and
$$
[t_1,\dots,t_l,t_{l+1},\dots,t_r]=
[1,\dots,1,\chi_1^{-1}(t')t_{l+1},\dots\chi_l^{-1}(t')t_r]
\text{ in } S,
$$
for $t'=(t_1,\dots,t_l)$. Now is clear that any point $[s]\in W(x_0)$ 
is fixed because $\bar O_{x_0}=\{z^\rho=0\mid\rho=l+1,\dots,n\}$, so 
that $W(x_0)\subset W^S$. But we could have described the action of 
$T$ on ${\mbb C}^r$ using the coordinates furnished by any other 
$n$-dimensional cone of $\Si$ and the conclusion would have been 
the same. 

\begin{proposition}{\label{prop:fix}}
The fixed point set of the $S$-action on $W$ is 
$$
W^S=\bigcup_{x\in X^S}W(x),\text{ with }
W(x):=\Gamma (C,\eP_\ud\times_\Cr\bar O_x)^o/T.
$$
Moreover, for any $x\in X^S$, $W(x)$ is defined by the fibre product
$$
\xymatrix{
W(x)\ar[d]\ar[r]&
\prod_{\rho\not\in\si_n(x)}{\mbb P}(\crl W_\rho)
\ar[d]\\ 
\eJ^r \ar[r]^-{\pr_x}&
\prod_{\rho\not\in\si_n(x)}\eJ.
}
$$
\end{proposition}

Before proceeding we notice that the $W(x)$'s are smooth and disjoint 
subvarieties of $W$, which is in agreement with the general result 
obtained in \cite{iv}.

\begin{proof}
The first part of the proposition being already proved, we are left 
with the second claim. Again, we are going to check it only for 
$x_0\in X^S$, the general statement coming from the symmetry of the 
problem. We observe that $Z_X\cap \bar O_{x_0}=\bar O_{x_0}\sm O_{x_0}
=\bigcup_\l (\{z^\l=0\}\cap\bar O_{x_0})$, because $O_{x_0}$ is the 
locus where $T$ acts freely. Therefore
$$
\begin{array}{rl}
\Gamma(C,\eP_\ud\times_\Cr\bar O_{x_0})^o\kern-1.5ex&
=\Gamma(C,\eP_\ud\times_\Cr\bar O_{x_0})\sm 
\Gamma(C,\eP_\ud\times_\Cr(Z_X\cap\bar O_{x_0}))\\[1.5ex] 
&
=\prod_{\l=1}^l\crl W_\l\sm 
\prod_{\l=1}^l\{s\mid \pr_{\crl W_\l}s=0\},
\end{array} 
$$ 
and the statement follows because $T\cong {({\mbb C}^*)}^l$ acts 
componentwise.
\end{proof}

From the proposition we see that no matter what $W$ looks like, 
its fixed point set for the torus action has a very down-to-earth 
description. A first byproduct is an explicit formula for the 
Euler number of the fibre $Y$.

\begin{corollary}{\label{Euler-Y}}
$$
\chi(Y)=\sum_{x\in X^S}
\prod_{\rho\not\in\si_n(x)}N_\rho.
$$
\end{corollary}

\begin{proof}
This equality is just a rewriting of the main result in \cite{iv}, 
but it can be proved in a more elementary way as follows: $Y$ being 
a toric variety, its Euler characteristic coincides with the number 
fixed points under the ${({\mbb C}^*)}^R/T$-action. Since this big 
torus contains $S$, the fixed point set is contained in the union of 
the $S$-fixed subvarieties. But fibrewise these are just products 
of projective spaces on which ${({\mbb C}^*)}^R$ acts in standard 
fashion.
\end{proof}

For applying the localization formula we must know the action of 
$S$ on the normal bundles to the fixed subvarieties. 

\begin{lemma}{\label{lm:exact}}
For any $x\in X^S$, the normal bundle of the fixed component $W(x)$ 
of $W$ fits in the following diagram:
\begin{align}{\label{diag-comm}}
\xymatrix{
 & & 0\ar[d] & 0\ar[d] & \\
0\ar[r]& 
{\cal O}(\Lie T)
\ar@{=}[d]\ar[r]& 
\bigoplus_{\rho\not\in\si_n(x)} q^*\crl W_\rho\otimes\L_\rho
\ar[d]\ar[r]&
\T_{W(x)/\eJ^r}
\ar[d]\ar[r]& 
0\\ 
0\ar[r]& 
{\cal O}(\Lie T)
\ar[r]& 
\bigoplus_{\rho=1}^r q^*\crl W_\rho\otimes\L_\rho
\ar[d]\ar[r]&
\T_{W/\eJ^r}
\ar[d]\ar[r]& 
0\\ 
 & & 
\bigoplus_{\rho\in\si_n(x)} q^*\crl W_\rho\otimes\L_\rho
\ar[d]\ar[r]^-\cong& 
\N_x:=\N_{W(x)|W}\ar[d] & \\
 & & 0 & 0 &
}
\end{align}
For the trivial action of $S$ on $\crl W_\rho$ and for the action
\begin{align}{\label{act-L}}
\begin{array}{l}
S\times\L_\rho\lar\L_\rho\text{ given by}\\[1.5ex] 
[t]\times[s,a]:=[t\times s,t_\rho a],\quad\forall\,
[t]\in S\text{ and }[s,a]\in\L_\rho,
\end{array}
\end{align}
all the homomorphisms in the diagram above are $S$-equivariant.
\end{lemma}

\begin{proof}
Since $\crl W^o\rar W$ is a principal $T$-bundle, we have 
the following $S$-equivariant exact sequence on $W$
\begin{align}{\label{t-inv}}
0\lar {\cal O}(\Lie T)\lar \T_{\crl W^o/\eJ^r}^{\rm\,inv}
\lar \T_{W/\eJ^r} \lar 0,
\end{align}
where $\T_{\crl W^o/\eJ^r}^{\rm\,inv}$ denotes the $S$-invariant 
relative tangent bundle to the total space of $\crl W^o$. But 
$\crl W^o$ is an open subset in a vector bundle over $\eJ^r$, so 
that the relative tangent bundle is canonically isomorphic to 
$Q^*\crl W=\oplus_\rho Q^*\crl W_\rho$, for $Q:\crl W\rar\eJ^r$ 
the projection. As $T$ preserves the decomposition of $\crl W$, 
$$
\T_{\crl W^o/\eJ^r}^{\rm\,inv}\cong 
Q^*\crl W/T=\bigoplus_{\rho=1}^r Q^*\crl W_\rho/T.
$$
We observe now that 
$Q^*\crl W_\rho/T\cong q^*\crl W_\rho\otimes\L_\rho$, the 
isomorphism being given by 
\begin{align}{\label{eq:iso}}
[s,w_\rho]\lmt w_\rho\otimes[s,1].
\end{align}
This proves the exactness of the middle row in the diagram 
\eqref{diag-comm}. A similar argument proves the exactness of 
the first horizontal sequence, and the last row is now a simple 
consequence. 

The very important thing which must be clarified yet is the way 
how $S$ acts on $q^*\crl W_\rho\otimes\L_\rho$. The sequence 
\eqref{t-inv} being $S$-equivariant, we have to  
describe the induced action on $q^*\crl W_\rho\otimes\L_\rho$ 
under the isomorphism \eqref{eq:iso}. For $[t]\in S$,
$$
[t]\times[s,w_\rho]=[t\times s,t_\rho w_\rho]
\lmt 
t_\rho w_\rho\otimes[t\times s,1]
=w_\rho\otimes[t\times s,t_\rho],
$$
so that we can see that indeed the $S$-action on $\crl W_\rho$ 
is trivial while the action on $\L_\rho$ is as in \eqref{act-L}.
\end{proof}

The next step is the computation of the equivariant first Chern 
classes for the restrictions of $\L_\rho$ to the fixed components 
$W(x)$. Before proceeding we notice that since 
$S=M^\vee\otimes_{\mbb Z}{\mbb C}^*$, there is a natural ring 
isomorphism $H^*(BS)\cong {\rm Sym}^\bullet M$, where $BS$ denotes 
as usual the classifying space for $S$.

\begin{lemma}{\label{lm:equiv-c-class}}
For $x\in X^S$, denote ${(u_\rho(x))}_{\rho\in\si_n(x)}\subset M$ 
the dual basis to ${(e^\rho)}_{\rho\in\si_n(x)}\subset M^\vee$ 
formed by the integral generators of $\si_n(x)$. Then 
$$
\begin{array}{l}
c_1^S(\L_\rho|_{W(x)})=\L_\rho|_{W(x)}+u_\rho(x),
\quad\forall\,\rho\in\si_n(x),\\[1ex]
c_1^S(\L_\rho|_{W(x)})=\L_\rho|_{W(x)},
\quad\forall\,\rho\not\in\si_n(x).
\end{array}
$$
\end{lemma}

\begin{proof}
It is clear that $c_1^S(\L_\rho|_{W(x)})=\L_\rho|_{W(x)}+u$, 
for some $u\in H^*(BS)$, and this element is precisely the 
weight of the action of $S$ on the stalk $\L_\rho|_{[s]}$ at 
some point $[s]\in W(x)$. Again, we shall make the computations 
for $x_0$ only: in this case the assignment 
$(\tau_{l+1},\dots,\tau_r)\mt [1,\dots,1,\tau_{l+1},\dots,\tau_r]$ 
gives an isomorphism ${({\mbb C}^*)}^n\xrar{\cong}S$ and is easy to 
see that
$$
{({\mbb C}^*)}^n\text{ acts on }\L_\rho\,
\left\{
\begin{array}{l}
\kern-.7ex\text{trivially, for }\rho=1,\dots,l
\text{ \it i.e. }\rho\not\in\si_n(x_0),\\[1ex]
\kern-.7ex\text{by }\tau_\rho,\text{ for }\rho=l+1,\dots,r
\text{ \it i.e. }\rho\in\si_n(x_0).
\end{array}\right.
$$
A short computation shows that the above isomorphism is induced 
precisely by the choice of the dual basis to $(e^{l+1},\dots,e^r)$, 
and the conclusion follows.
\end{proof}

We are finally in position to apply the localization formula 
for computing intersection numbers. For positive integers 
$m_1,\dots,m_r$, we want is to compute the push-forward 
$q_*\bigl(\L_1^{m_1}\cdot{\dots}\cdot\L_r^{m_r}\bigr)
\in H^*(\eJ^r)$; 
for shorthand, we write $\phi$ for this product so that we must 
compute $q_*\phi$. The case of interest for us is when 
$$
m_1+\dots+m_r=\dim V=N_1+\dots+N_r+l(g-1).
$$
We denote $\tld e(\N_x)$ the equivariant Euler 
characteristic of $\N_x\rar W(x)$ and let $\tld\phi$ to be the 
equivariant analog of $\phi$. If $\jm_x:W(x)\hra W$ is the 
inclusion, the localization formula for $\tld\phi$ reads
$$
\tld\phi
=\sum_{x\in X^S}
{(\jm_x)}_*\frac{\jm_x^*\tld\phi}{\tld e(\N_x)}.
$$
Since $S$ acts trivially on $\eJ^r$, composing with the projection 
$q:W\rar\eJ^r$ we obtain
\begin{align}{\label{eq:push-forward}}
q_*\tld\phi
=\sum_{x\in X^S}
{(q_x)}_*\frac{\jm_x^*\tld\phi}{\tld e(\N_x)}
\in H^*_S(\eJ^r)=H^*(BS)\otimes H^*(\eJ^r),
\end{align}
for $q_x:=q\circ\jm_x$, so that $q_*\phi$ is the $H^*(BS)$-free 
term in $q_*\tld\phi$. 
Lemmas \ref{lm:exact} and \ref{lm:equiv-c-class} imply that 
$$\kern-.5ex
\begin{array}{rl}
\tld e(\N_x)\kern-1.5ex&
\disp
=
\kern-1.5ex
\prod_{\rho\in\si_n(x)}\kern-.8ex
\biggl( 
\kern-.3ex
(u_\rho(x)
\kern-.1ex+\kern-.1ex
\L_\rho)^{N_\rho}
\kern-.3ex
-
\th_\rho 
(u_\rho(x)
\kern-.1ex+\kern-.1ex
\L_\rho)^{N_\rho-1}
\kern-.3ex
+
\frac{\th_\rho^2}{2!}
(u_\rho(x)
\kern-.1ex+\kern-.1ex
\L_\rho)^{N_\rho-2}
\kern-.3ex
-
\dots
\kern-.3ex
\biggr)\\[1.5ex]
&
\disp\kern-1.5ex
\srel{_{N_\rho > g}}{=}
\kern-1.5ex
\prod_{\rho\in\si_n(x)}
(u_\rho(x)+\L_\rho)^{N_\rho}
\cdot
\exp\biggl(-\frac{\th_\rho}{u_\rho(x)+\L_\rho}\biggr),
\end{array}
$$
where we are using the fact (see for instance page 336 
in \cite{acgh}) that the total Chern class of the Picard bundle 
$\crl W_\rho\rar\eJ$ is $c(\crl W_\rho)=\exp(-\th_\rho)$, with 
$\th_\rho$ the class of the theta divisor (the lower index $\rho$ 
indicates that we are on the $\rho^{\rm th}$ copy of $\eJ$ in 
$\eJ^r$). As a consequence,
\begin{align}{\label{eq:jx}}
\begin{array}{rl}
\disp
\frac{\jm_x^*\tld\phi}{\tld e(\N_x)}\kern-1.5ex
&
\disp
=\frac{\disp
\prod_{\rho\not\in\si_n(x)}\L_\rho^{m_\rho}
\cdot
\prod_{\rho\in\si_n(x)}(u_\rho(x)+\L_\rho)^{m_\rho}
}
{\disp
\prod_{\rho\in\si_n(x)}
(u_\rho(x)+\L_\rho)^{N_\rho}
\cdot
\exp\biggl(-\frac{\th_\rho}{u_\rho(x)+\L_\rho}\biggr)
}
\\[2ex]
&
\disp
=\prod_{\rho\not\in\si_n(x)}\L_\rho^{m_\rho}
\cdot
\prod_{\rho\in\si_n(x)}
(u_\rho(x)+\L_\rho)^{m_\rho-N_\rho}
\cdot
\exp\biggl(\frac{\th_\rho}{u_\rho(x)+\L_\rho}\biggr).
\end{array}
\end{align}
Plugging this into equality \eqref{eq:push-forward} we find 
the formula for the push forward of the class $\phi$.

\begin{theorem}{\label{thm:push-f}}
For any positive integers $m_\rho$, 
$q_*\bigl(
\L_1^{m_1}\cdot{\dots}\cdot\L_r^{m_r}
\bigl)$ 
is the constant term of the polynomial 
\begin{align}{\label{sum}}
\sum_{x\in X^S}
{(q_x)}_*
\left[
\prod_{\rho\not\in\si_n(x)}
\L_\rho^{m_\rho}
\cdot\kern-1ex
\prod_{\rho\in\si_n(x)}
(u_\rho(x)+\L_\rho)^{m_\rho-N_\rho}
\cdot
\exp\biggl(\frac{\th_\rho}{u_\rho(x)+\L_\rho}\biggr)
\right].
\end{align}
\end{theorem}

It is certainly not apparent that this sum is a polynomial 
expression in the formal variables $u_\rho(x_0)=:u_\rho$, 
and computing its constant term is not an easy task 
(we must chose a basis of $M$ for writing \eqref{sum} as an 
element of ${\rm Sym}^\bullet M$, and $x_0\in X$ is already 
our favorite fixed point in $X$). In most cases it is not 
true that the constant term of the whole sum is the sum of 
the individual constant terms. The question we are going 
to discuss is how to apply the result in concrete situations? 
Expanding the terms in \eqref{sum} is a straightforward 
computation, but requires a little patience.

\begin{lemma}{\label{lm:patience}}
For $p\geq 0$,
$$\kern-1ex
\begin{array}{ll}\disp 
{(u+\L)}^p\exp\biggl(\frac{\th}{u+\L}\biggr)\kern-1.5ex
&\disp
=\sum_{0\leq k\leq p}
\left[
\sum_{0\leq a\leq k}\binom{p-a}{k-a}\L^{k-a}\frac{\th^a}{a!}
\right]
u^{p-k}\\[4ex]
&\disp
+
\th^{p+1}
\kern-.25ex
\sum_{0\leq k}
\kern-.25ex
\left[
\sum_{0\leq b\leq k}{(-1)}^{k-b}\binom{k}{k-b}\L^{k-b}
\frac{\th^b}{(p+1+b)!}
\right]
\kern-.5ex
\frac{1}{u^{k+1}},
\end{array}
$$
while for $p\geq 1$,
$$
\frac{1}{{(u+\L)}^p}\exp\biggl(\frac{\th}{u+\L}\biggr)
=
\sum_{0\leq k}
\left[
\sum_{0\leq b\leq k}
{(-1)}^{k-b}\binom{p+k-1}{k-b}\L^{k-b}\frac{\th^b}{b!}
\right]
\frac{1}{u^{p+k}}.
$$
\end{lemma}

The lemma implies that the expression \eqref{sum} is a sum of 
rational functions, which are quotients of homogeneous polynomials 
in  the $u_\rho(x)$'s. Its constant term is obtained by adding the 
functions having total degree zero (the total degree is the 
difference between the degrees of the nominator and the denominator). 

For every $x\in X^S$, the coefficients of these functions are 
products in $\L_\rho|_{W(x)}$, and their push-forward by $q_x$ 
can be easily computed because, for $\rho\not\in\si_n(x)$, the 
restriction of $\L_\rho|_{W(x)}$ is simply (the pull-back of) the 
usual relatively ample line bundle 
${\cal O}_{{\mbb P}(\crl W_\rho)}(1)\rar{\mbb P}(\crl W_\rho)$ 
while the remaining line bundles $\L_\rho|_{W(x)}$, 
$\rho\in\si_n(x)$, are linear combinations of the previous ones 
(see equation \eqref{lin-rel}). 

The integrals which appear are of the type
$
{(q_x)}_*\bigl(
\prod_{\rho\not\in\si_n(x)}\L_\rho^{k_\rho}
\bigr),
$
and is quite known that
\begin{align}{\label{segre}}
{(q_x)}_*\biggl(
\prod_{\rho\not\in\si_n(x)}\L_\rho^{k_\rho}
\biggr)
=
\begin{cases}
\disp\prod_{\rho\not\in\si_n(x)}
\frac{\th_\rho^{k_\rho-N_\rho+1}}{(k_\rho-N_\rho+1)!}
&\text{if $N_\rho-1\leq k_\rho\leq N_\rho+g-1\,\forall\rho$,}\\ 
\kern6em 0&\text{otherwise}.
\end{cases}
\end{align}
Now we recall that we are actually interested in computing 
integrals on $V$, the compactification of the space 
$\Mor_\ud(C,X)$, so that we must pull-back to $\eJ^l$ the class 
$q_*\phi\in H^*(\eJ^r)$, using the morphism $\psi$ defined by 
\eqref{map-phi}. But $\psi$ is explicitely given in terms of 
the combinatorics of the fan $\Si$. 

The difficulty in applying formula \eqref{sum} relies in the 
fact that we are not allowed to set to zero the variables 
${\{ u_\rho(x)\}}_{\rho,x}$. However, as we will see in a moment, 
for special choices of the exponents ${(m_\rho)}_\rho$ this is 
possible, and in this cases we obtain very explicit formulae for 
the corresponding intersection products. 

Let us consider positive integers $a_1,\dots,a_r$ with the 
property that $a_1+{\dots}+a_{r-1}-a_r=0$ and take 
\begin{align}{\label{eq:choice}}
\begin{array}{l}
m_\rho=N_\rho+g+a_\rho=d_\rho+1+a_\rho,
\text{ for } \rho=1,\dots,r-1,\text{ and}
\\[1ex] 
m_r=N_r-(n-1)g-l-a_r=d_r-ng-(l+a_r-1).
\end{array}
\end{align}
Of course, such a choice is possible only when $d_r$ is 
large enough for $m_r$ to be positive. For such a choice, 
the total degree of the functions which appear in the products 
corresponding to the fixed points $x\in X^S$ with $r\in\si_n(x)$ 
is strictly negative; therefore these terms do not contribute to 
the intersection product. On the other hand, the products 
corresponding to $x\in X^S$ with $r\not\in\si_n(x)$ are honest 
polynomials, so that we are allowed to set the variables to zero. 

\begin{proposition}{\label{prop:explicit}}
For integers ${(m_\rho)}_\rho$ as in \eqref{eq:choice}, 
$$
q_*(\L_1^{m_1}\cdot{\dots}\cdot\L_r^{m_r})
=
\sum_{x\in X^S\!,\, r\not\in\si_n(x)}
\kern-1.5ex
{(q_x)}_*
\left[
\prod_{\rho\not\si_n(x)}\L_\rho^{m_\rho}\cdot
\prod_{\rho\in\si_n(x)}
\biggl(
\sum_{b=0}^{m_\rho-N_\rho}
\L_\rho^{m_\rho-N_\rho-b}\frac{\th_\rho^b}{b!}
\biggr)
\right].
$$
\end{proposition}

Analogous formulae can be obtained for any other $k\in\{1,\dots,r\}$ 
or, when it is possible, by taking several $k$'s such that $m_k-N_k<0$ 
in a suitable way.

We conclude this section with a vanishing result, which is a direct 
consequence of theorem \ref{thm:push-f}.

\begin{proposition}{\label{prop:0}}
Consider $J\subset\Si(1)$ such that the vectors 
${(e^\rho)}_{\rho\in J}$ do not generate a cone of $\Si$ 
(the primitive collections are the smallest subsets of 
$\Si(1)$ with this property). If ${(m_\rho)}_\rho$ are 
positive integers such that $m_\rho\geq N_\rho+g$ for all 
$\rho\in J$, then 
$$
\int_V\L_1^{m_1}\cdot{\dots}\cdot\L_r^{m_r}=0.
$$
\end{proposition}

\begin{proof}
Indeed, for such a set $\bigcap_{\rho\in J}D_\rho=\emptyset$ and 
consequently $J\sm \si_n(x)\neq\emptyset$ for any $x\in X^S$ 
(otherwise $x\in\bigcap_{\rho\in J}D_\rho$, a contradiction). 
The conclusion follows now from the relations \eqref{segre}.
\end{proof}

We should point out that it is not always possible to chose 
integers with the property above: examples in this sense are 
the projective spaces. On the other hand, when there are 
primitive collections with less than $l={\rm rank\,}T$ 
elements, there are integers having the desired property.

The conclusion of this article is that we were able to translate 
the integration problem on the space of morphisms from the curve 
$C$ into the toric variety $X$ in an integration problem on a 
power of the Picard torus of $C$. Interestingly enough, these 
integrals depend only on the combinatorics of the fan defining 
$X$ and on the theta classes of the Picard varieties.


\end{document}